\documentclass[12pt]{amsart}

\title[Embedding relatively hyperbolic groups]{Embedding relatively hyperbolic groups in products of trees}
\author{John M. Mackay}
\address{Mathematical Institute \\
 University of Oxford \\ Oxford, UK.}
\email{john.mackay@maths.ox.ac.uk}
\author{Alessandro Sisto}
\email{sisto@maths.ox.ac.uk}
\date{March 20, 2013}

\usepackage{amsmath,amsthm,amsfonts,graphicx,amssymb}
\usepackage{array}
\usepackage{hyperref}

\numberwithin{equation}{section}
\newtheorem{theorem}[equation]{Theorem}
\newtheorem{proposition}[equation]{Proposition}
\newtheorem{corollary}[equation]{Corollary}
\newtheorem{lemma}[equation]{Lemma}

\newtheorem{definition}[equation]{Definition}
\newtheorem{remark}[equation]{Remark}

\newtheoremstyle{citing}
  {3pt}
  {3pt}
  {\itshape}
  {}
  {\bfseries}
  {}
  {.5em}
  {\thmnote{#3}}

\theoremstyle{citing}
\newtheorem*{varthm}{}

\DeclareMathOperator{\diam}{diam}

\DeclareMathOperator{\asdim}{asdim}

\newcommand{\lasdim}{{\ell\text{-}\!\asdim}}
\newcommand{\ldim}{{\ell\text{-}\!\dim}}
\DeclareMathOperator{\ecodim}{eco-dim}

\newcommand{\eps}{\epsilon}

\newcommand{\bdry}{\partial_\infty}

\newcommand{\cH}{\mathcal{H}}
\newcommand{\cT}{\mathcal{T}}

\newcommand{\cC}{\mathcal{C}}

\newcommand{\cP}{\mathcal{P}}
\newcommand{\cG}{\mathcal{G}}

\newcommand{\ra}{\rightarrow}
\newcommand{\R}{\mathbb{R}}

\newcommand{\N}{\mathbb{N}}
\newcommand{\Z}{\mathbb{Z}}
\newcommand{\HH}{\mathbb{H}}

\newcommand{\Nil}{\mathrm{Nil}}
\newcommand{\Sol}{\mathrm{Sol}}

\newcommand{\cU}{\mathcal{U}}

\def\XXint#1#2#3{{\setbox0=\hbox{$#1{#2#3}{\int}$}
\vcenter{\hbox{$#2#3$}}\kern-.5\wd0}}

\numberwithin{equation}{section}

\begin{document}

\begin{abstract}
	 We show that a relatively hyperbolic group quasi-isometrically embeds in a product of finitely many trees if the peripheral subgroups do, and we provide an estimate on the minimal number of trees needed. Applying our result to the case of 3-manifolds, we show that fundamental groups of closed 3-manifolds have linearly controlled asymptotic dimension at most 8. To complement this result, we observe that fundamental groups of Haken 3-manifolds with non-empty boundary have asymptotic dimension~2.
\end{abstract}

\maketitle

\section{Introduction}\label{sec-intro}

If a group admits a quasi-isometric embedding into the product of
finitely many trees then one can say a lot about the large scale
geometry of the group.  For example, this bounds the (linearly
controlled) asymptotic dimension of the group.  It has implications for
the Hilbert compression of the group, and also the topological dimension
of any asymptotic cone.  Such embeddings have recently been constructed
for hyperbolic groups \cite{Bu-05-asdim-capdim, BuDrSc-07-hyp-grp-prod-trees}.

In this paper we show that if every peripheral group of a relatively
hyperbolic group embeds into the product of finitely many trees, then
the entire group does.  Moreover, we give explicit bounds for the number
of trees required.

We begin by recalling what is known for hyperbolic groups.
In this case, the situation has been completely understood by work of Buyalo and Lebedeva
\cite[Theorems 1.5, 1.6]{BuLe-selfsim}.  (For further refinements, see \cite{BuDrSc-07-hyp-grp-prod-trees}.)
\begin{theorem}\label{thm-buyalo-schroeder}
	Let $G$ be a Gromov hyperbolic group.  Then $G$ admits a quasi-isometric embedding into the product
	of $n+1$ metric trees, where $n = \dim \bdry G$ is the topological dimension of the boundary.
	Moreover, $G$ does not embed into any product of $n$ metric trees.
\end{theorem}

Our main result is the following.  
(For the definition of a relatively hyperbolic group, see Section~\ref{sec-rel-hyp-defs}.)
\begin{theorem}\label{thm-main1}
	Suppose the group $G$ is hyperbolic relative to subgroups $H_1, \ldots, H_n$.
	If each $H_i$ quasi-isometrically embeds into a product of $m$ metric trees, then
	$\asdim(G) < \infty$ and 
	$G$ quasi-isometrically embeds into a product of $M$ metric trees, where
	\[
		M = \max\{\asdim(G), m+1\}+m+1 < \infty.
	\]
	Conversely, if $G$ quasi-isometrically embeds into a product of $N$ metric trees, then each peripheral group does also.
\end{theorem}
The assertion $\asdim(G) < \infty$ follows from \cite{Osin-asdimrelhyp}.

Both these theorems use results that study the relationship between the asymptotic dimension of a hyperbolic space
and the linearly controlled metric dimension of its boundary.  We now proceed to define these and related concepts.
(See \cite{BuSc-07-asymp-geom-book} for more discussion.)

Suppose $\cU$ is a family of subsets of a metric space $X$.
We say $\cU$ is \emph{$D$-bounded} if the diameter of every $U \in \cU$ is at most $D$.
The \emph{$r$-multiplicity} of $\cU$ is the infimal integer $n$ so that every subset of $X$ with diameter
less than or equal to $r$ meets at most $n$ subsets of $\cU$.

The \emph{asymptotic dimension} of $X$, denoted by $\asdim(X)$, is the smallest $n \in \N \cup \{\infty\}$ so that
for all $r > 0$, there exists $D(r) < \infty$ and a $D(r)$-bounded cover $\cU$ of $X$ with $r$-multiplicity at most $n+1$
\cite[1.E]{Gro-91-asymp-inv}.
The \emph{linearly controlled asymptotic dimension} of $X$, denoted by $\lasdim(X)$, is the smallest $n$
so that there exists $L < \infty$ with the property that for all sufficiently large $r<\infty$, there exists an $Lr$-bounded 
cover $\cU$ of $X$ with $r$-multiplicity at most $n+1$ \cite{LaSc-05-nagata,BDHM-09-asdiman-via-lip}.
(This is sometimes called ``asymptotic Assouad-Nagata dimension'' in the literature.)

The following definition simplifies our discussion of how metric spaces embed in trees.
(For the related notion of ``$t$-rank'', see \cite[page 176]{BuSc-07-asymp-geom-book}.)
\begin{definition}\label{def-ecodim}
Given a metric space $X$, let $\ecodim(X)$ be the smallest $n \in \N$ so that 
$X$ quasi-isometrically embeds in the product of $n$ metric trees, and
set $\ecodim(X)= \infty$ if no such embedding exists.
\end{definition}
The inequalities $\asdim(X) \leq \lasdim(X) \leq \ecodim(X)$ hold for any metric space $X$. 
These three quantities are equal for hyperbolic groups, as shown in \cite{BuLe-selfsim}.
Both equalities can fail in general: there are groups $G$ with finite asymptotic dimension but infinite
linearly controlled asymptotic dimension~\cite{Nowak-07-asdim-an}.
The discrete Heisenberg group $H$ has $\asdim(H)=\lasdim(H)=3$,
but $\ecodim(H)=\infty$ (see the discussion in Section~\ref{sec-three-mfld}).

As an aside, Lang and Schlichenmaier show that if $\lasdim(X, d) < \infty$, then for sufficiently
small $\eps>0$, the snowflaked space $(X, d^\eps)$ has
$\ecodim(X, d^\eps) \leq \lasdim(X, d)+1$ \cite[Theorem 1.3]{LaSc-05-nagata}.

For hyperbolic groups, these asymptotic invariants are related to the local properties of the boundary.
The \emph{linearly controlled metric dimension} of $X$, 
denoted by $\ldim(X)$, is the smallest $n$ so that there exists $L < \infty$ with the property that for 
all sufficiently small $r>0$, there exists an $Lr$-bounded 
cover $\cU$ of $X$ with $r$-multiplicity at most $n+1$ \cite[Prop.\ 3.2]{Bu-05-asdim-capdim}.  (This is also
referred to as ``capacity dimension''.)
Buyalo shows the following embedding theorem.
\begin{theorem}[{\cite[Theorem 1.1]{Bu-05-cdim-prod-trees}}]\label{thm-buyalo-cdim-trees}
	Suppose $X$ is a visual, Gromov hyperbolic metric space, and $\ldim(\bdry X) < \infty$.
	Then $\ecodim(X) \leq \ldim(\bdry X) +1$.
\end{theorem}

The following proposition gives a simple bound for the linearly controlled metric dimension of the boundary of $X$.
\begin{varthm}[Proposition \ref{prop-cdim-asdim}]
	Suppose $X$ is a Gromov hyperbolic geodesic metric space.  Then $\ldim(\bdry X) \leq \asdim(X)$.
\end{varthm}
This proposition does not seem to be recorded in the literature, possibly because in the case when the space
admits a cocompact isometric action the stronger equality $\asdim(X) = \ldim(\bdry X) + 1$ holds \cite{BuLe-selfsim}.
Even without such an action, Buyalo shows the inequality $\asdim(X) \leq \ldim(\bdry X) +1$ \cite{Bu-05-asdim-capdim}.

In the case of a relatively hyperbolic group $(G,\{H_i\})$, a natural choice for $X$ is $X(G)$, the Bowditch space associated to
$(G,\{H_i\})$, see Definition~\ref{def-rel-hyp}.  This is a visual, Gromov hyperbolic, geodesic metric space.
We bound the asymptotic dimension of $X(G)$ in terms of the asymptotic dimension
of $G$ and the linearly controlled asymptotic dimension of the peripheral groups.
\begin{varthm}[Proposition \ref{prop-asdim-bowditch}]
	Let $G$ be hyperbolic relative to $H_1,\dots, H_n$, and
	let $m = \max_{i=1, \ldots, n} \lasdim(H_i)$.
	Then
	\[
		\max\{\asdim(G),m\} \leq \asdim(X(G)) \leq \max \{\asdim(G), m+1\}.
	\]
\end{varthm}
Osin had earlier shown that $\asdim(G)$ is finite if $\asdim(H_i) < \infty$ for each $i$ \cite[Theorem 1.2]{Osin-asdimrelhyp}.

The product of $n$ (unbounded) metric trees has linearly controlled asymptotic dimension $n$ (see, e.g., \cite{LaSc-05-nagata}).
Therefore, Theorem~\ref{thm-buyalo-cdim-trees} and Propositions~\ref{prop-cdim-asdim} and \ref{prop-asdim-bowditch}
combine to show the following.
\begin{corollary}\label{cor-bowditch-trees}
	Suppose the group $G$ is hyperbolic relative to $H_1, \ldots, H_n$.
	Let $m = \max \{ \ecodim(H_i) \}$.
	Then $\ecodim(X(G)) \leq \max\{\asdim(G), m+1\}+1$.
\end{corollary}

In order to prove Theorem~\ref{thm-main1}, we use work of Bestvina, Bromberg and Fujiwara to combine the
embeddings of the peripheral groups and of $X(G)$ into a single embedding of $G$ into a product of trees.
\begin{varthm}[Theorem \ref{prodqtrees}]
	 Let $G$ be hyperbolic relative to $H_1,\dots, H_n$, and suppose that each $H_i$ admits a quasi-isometric 
	 embedding into the product of $m$ trees. Then $G$ admits a quasi-isometric embedding into the product of
	 $m$ trees and either $X(G)$ or the coned-off graph $\hat{G}$.
\end{varthm}
This theorem completes the proof of Theorem~\ref{thm-main1}.  Note that the converse statement is automatic,
as peripheral groups of a relatively hyperbolic group are undistorted in the ambient group
\cite[Lemma 4.15]{DSp-05-asymp-cones}.

Rather than using the Bowditch space $X(G)$, to improve the bounds in Theorem~\ref{thm-main1}
one might hope to use instead the (hyperbolic) coned-off graph $\hat{G}$.  
However, the non-locally finite nature of $\hat{G}$ leads to a non-compact boundary, and
much of the machinery used above no longer applies.

The embeddings we consider are not and cannot required to be equivariant, 
because there exist hyperbolic groups with property (T) (for example, any cocompact lattice in $Sp(n,1)$),
which in particular cannot act interestingly on trees. 
We remark that the inverse problem of equivariantly embedding trees into hyperbolic spaces is 
treated in~\cite{BIM-05-embed-trees-in-hyp}.

Theorem~\ref{thm-main1} has the following immediate corollary.
\begin{corollary}
	If the peripheral groups of a relatively hyperbolic group $G$ each quasi-isometrically embed 
	into the product of finitely many trees, then $\lasdim(G) < \infty$.
\end{corollary}
If the group $G$ satisfies $\lasdim(G) < \infty$ then it has Hilbert compression $1$ \cite[Theorem 1.1.1]{Gal-08-asdim-hilbert-compression}.
For more general estimates on the Hilbert compression of a relatively hyperbolic group, see \cite{Hume}.

In Section~\ref{sec-three-mfld} we apply our results to $3$-manifold groups, and show, amongst other results,
the following.
\begin{varthm}[Theorem \ref{thm-three-mfld}]
	Let $G = \pi_1(M)$, where $M$ is a compact, orientable $3$-manifold whose (possibly empty) boundary is a union of tori.
	Then $\ecodim(G) < \infty$ if and only if no manifold in the prime decomposition of $M$ has $\mathrm{Nil}$ geometry;
	in this case, $\ecodim(G) \leq 8$.
	
	In any case, $\lasdim(G) \leq 8$.
\end{varthm}

Finally, in Appendix \ref{sec-hnn} we prove a bound on the asymptotic dimension of HNN extensions,
following work of Dranishnikov.

\subsection{Recent developments}
Since the appearance of the first version of this paper, Hume \cite{Hu-bbf-qtrgr} has shown that the construction of Bestvina, Bromberg and Fujiwara that we use in Section \ref{sec-embed-prod-trees-and-bow} gives rise to spaces that are quasi-isometric to a tree-graded space. Tree-graded spaces are defined in \cite{DSp-05-asymp-cones}; roughly speaking, they are tree-like arrangements of certain specified subspaces called pieces (for example, consider a Cayley graph of a free product). It is known \cite{BrHi-09-AN-tree-graded} that if $X$ is tree-graded and $\cP$ is its collection of pieces then
$$\lasdim(X)\leq\max_{P\in\cP}\{1,\lasdim(P)\}$$
when, say, there are finitely many $(L,C)$-quasi-isometry types of pieces, for some given $L,C$.

As noticed in \cite{Hu-bbf-qtrgr}, the arguments in Section \ref{sec-embed-prod-trees-and-bow} (just considering $\pi_Y$ instead of $f_{i,Y}\circ \pi_Y$ in Subsection \ref{ssec-prodqtrees}), combined with \cite{Hu-bbf-qtrgr}, give the following.

\begin{theorem}[cf.\ Theorem~\ref{prodqtrees}]
 Let $G$ be hyperbolic relative to $H_1,\dots,H_n$. Then $G$ quasi-isometrically embeds in the product of a space $\cT$ and either the coned-off graph $\hat{G}$ or the Bowditch space $X(G)$, where $\cT$ is tree-graded and each piece is uniformly quasi-isometric to some $H_i$.
\end{theorem}

In particular, once again as noticed in \cite{Hu-bbf-qtrgr}, using the aforementioned result from \cite{BrHi-09-AN-tree-graded} one can get the following.

\begin{corollary}
 A relatively hyperbolic group has finite $\lasdim$ if and only if its peripheral subgroups do.
\end{corollary}

\subsection*{Acknowledgements} 
We thank Alexander Dranishnikov, David Hume, Urs Lang, Enrico Le Donne and the referee for helpful comments, 
and Marc Lackenby for providing a fix for the proof of Proposition \ref{nonclosed}.

\section{Relatively hyperbolic groups}\label{sec-rel-hyp-defs}

In this section we define relatively hyperbolic groups and their (Bow\-ditch) boundaries.

\subsection{Definitions}
 There are many definitions of relatively hyperbolic groups.
 We will give one in terms of actions on a cusped space.
 First we define our model of a horoball.
 
 \begin{definition}
 	Suppose $\Gamma$ is a connected graph with vertex set $V$ and edge set $E$,
 	where every edge has length one.
 	The \emph{horoball} $\cH(\Gamma)$ is defined to be the graph with
 	vertex set $V \times \N$ and
 	edges $((v, n), (v, n+1))$ of length $1$, for all $v \in V$, $n \in \N$,
 	and edges $((v, n), (v', n))$ of length $e^{-n}$, for all $(v, v') \in E$.	
 \end{definition}
 
 Note that $\cH(\Gamma)$ is quasi-isometric to the metric space
 constructed from $\Gamma$ by gluing to each edge in $E$ a copy of the strip
 $[0,1] \times [1, \infty)$ in the upper half-plane model of $\HH^2$,
 where the strips are attached to each other along $v \times [1,\infty)$.
 
 As is well known, these horoballs are hyperbolic with boundary a single point.
 Moreover, it is easy to estimate distances in horoballs. We will write $A\approx B$ if the quantities $A$ and $B$ differ by some constant.
 \begin{lemma}\label{distestim}
 	Suppose $\Gamma$ and $\cH(\Gamma)$ are defined as above.
 	Let $d_\Gamma$ and $d_\cH$ denote the corresponding path metrics.
 	Then for each $(x,m), (y,n) \in \cH(\Gamma)$, we have 
 	\[
 		d_\cH((x,m),(y,n)) \approx 2 \ln(d_\Gamma(x,y)e^{-\min\{m,n\}}+1)+|m-n|.
 	\]
 \end{lemma}
 \begin{proof}
	We may assume that $m \leq n$.
 	We can suppose $d_\Gamma(x,y)\geq e^m$, as if not $|m-n| \leq d_\cH((x,m),(y,n))\leq 1+|m-n|$. In particular
 	$\ln(d_\Gamma(x,y)/e^m+1)\approx \ln d_\Gamma(x,y)-m$. By construction of $\cH(\Gamma)$, any geodesic 
 	$\gamma$ between $(x,m)$ and $(y,n)$ in $\cH(\Gamma)$ must go from $(x,m)$ to $(x,t)$ changing
 	only the second coordinate,
 	then follow a geodesic $\gamma' \subset \Gamma \times \{t\} \subset \cH(\Gamma)$
 	to $(y,t)$, then back to $(y,n)$.	
 	Thus
 	\begin{equation}\label{eq-horoball-distance-estimate}
 		d_\cH((x,m),(y,n)) \leq 2(t-m)+|n-m|+ e^{-t} d_\Gamma(x,y),
 	\end{equation}
 	with equality for the best choice of $t$. It is readily seen that this value is attained for the least $t$ so that $l_t=e^{-t} d_\Gamma(x,y)$ satisfies $l_t/e+2\geq l_t$, that is, $l_t\leq 2e/(e-1)=\rho$.
 	So the best choice of $t$ is $t = \lceil \ln (d_\Gamma(x,y)/\rho)\rceil$, 
 	and the right hand side of \eqref{eq-horoball-distance-estimate}
 	is $2(\ln d_\Gamma(x,y)-m)-2\ln\rho + |n-m| +\epsilon$, where $|\epsilon|\leq 2+\rho$.
 \end{proof}

\begin{definition}\label{def-rel-hyp}
	Suppose $G$ is a finitely generated group, and $\{H_i\}_{i=1}^n$
	a collection of finitely generated subgroups of $G$.
	Let $S$ be a finite generating set for $G$, so that $S \cap H_i$ generates
	$H_i$ for each $i=1, \ldots, n$.
	
	Let $\Gamma(G,S)$ be the Cayley graph of $G$ with respect to $S$.
	Let $X(G) =  X(G, \{H_i\}, S)$ be the space resulting from gluing to $\Gamma(G,S)$
	a copy of $\cH(\Gamma(H_i, S \cap H_i))$ to each coset $g H_i$ of
	$H_i$, for each $i= 1, \ldots, n$.  We call $X(G)$ the \emph{Bowditch space} associated to $(G, \{H_i\}, S)$.
	
	We say that $(G, \{H_i\})$ is \emph{relatively hyperbolic} if $X(G)$ is Gromov hyperbolic,
	and call the members of $\{H_i\}$ \emph{peripheral subgroups}.
\end{definition}

This is equivalent to the other usual definitions of (strong) relative hyperbolicity;
see \cite{Bow-99-rel-hyp}, \cite[Theorem 3.25]{Gro-Man-08-dehn-rel-hyp}.

\subsection{Visual metric}
\label{visual-metric}
Let $X$ be a geodesic, Gromov hyperbolic space (not necessarily proper), with fixed base point $0 \in X$.
Suppose all geodesic triangles are $\delta$-slim. 
One equivalent definition of the boundary $\bdry X$ is as the set of equivalence classes of $(1,20\delta)$-quasigeodesic 
rays $\gamma: [0, \infty) \ra X$, with $\gamma(0)=0$, where two rays are equivalent if they are at finite Hausdorff distance
from each other.
Let $(x | y) = (x | y)_0$ denote the Gromov product on $\bdry X$ with respect to $0$.
Up to an additive error, $(x|y)$ equals the distance from $0$ to
some (any) $(1,20\delta)$-quasigeodesic line from $x$ to $y$ \cite[Remark 2.16]{KB-02-boundaries-survey}.

A metric $\rho$ on $\bdry X$ is a visual metric if 
there exists $C_0, \eps >0$ so that $\frac{1}{C_0} e^{-\eps (x|y)} \leq \rho(x,y) \leq C_0 e^{-\eps (x|y)}$ 
for all $x, y \in \bdry X$. Boundaries of Bowditch spaces will always be endowed with a visual metric.
The linearly controlled metric dimension of such boundary is independent of the choice of visual metric \cite{Bu-05-cdim-prod-trees}.

For more results on the geometric properties of $\bdry X(G)$, see \cite{MS-12-relhypplane}.

\subsection{Distance formula}
Let $G$ be a relatively hyperbolic group and let ${\bf Y}$ be the collection of all left cosets of peripheral subgroups. 
Fix a Cayley graph of $G$ with path metric $d$.
For $Y\in{\bf Y}$, let $\pi_Y:G \ra Y$ be a closest point projection map onto $Y$ with respect to $d$.
Denote by $\hat{G}$ the \emph{coned-off graph} of $G$, that is to say the metric graph obtained from the Cayley graph of $G$ by adding an edge 
of length one connecting each pair of (distinct) vertices contained in the same left coset of peripheral subgroup.
The path metric on $\hat{G}$ is denoted by $d_{\hat{G}}$.
Let $\big\{\big\{x\big\}\big\}_L$ denote $x$ if $x>L$, and $0$ otherwise.
We write $A\approx_{\lambda,\mu} B$ if $A/\lambda-\mu\leq B\leq \lambda A+\mu$. The following is proven in \cite{projrelhyp}.
\begin{theorem}[Distance formula for relatively hyperbolic groups]\label{distform}
There exists $L_0$ so that for each $L\geq L_0$ there exist $\lambda,\mu$ so that the following holds. If $x,y\in G$ then
\begin{equation}\label{eqn-distform}
	d(x,y)\approx_{\lambda,\mu} \sum_{Y\in\bf{Y}} \big\{\big\{d(\pi_Y(x),\pi_Y(y))\big\}\big\}_L+d_{\hat{G}}(x,y).
\end{equation}
\end{theorem}

\section{Large scale and boundary dimension estimates}\label{sec-dim-estimates}

In this section we bound the asymptotic dimension of the Bowditch space of a relatively hyperbolic group.
We also bound the linearly controlled metric dimension of the boundary of a Gromov hyperbolic space by its asymptotic dimension.

Observe that at the cost of a slight relaxation in the value of $r$,
a collection of subsets $\cU$ has $r$-multiplicity at most $m$ if and only if we can write
$\cU = \bigcup_{i=1}^m \cU_i$, where each $\cU_i$ is \emph{$r$-disjoint}: it is a collection of subsets pairwise separated by 
a distance of at least $r$.
This follows from an application of Zorn's lemma.
We will use this alternative characterization throughout this section.

\subsection{Asymptotic dimension of Bowditch spaces}

First we bound the asymptotic dimension of a horoball.
\begin{proposition}\label{asdimhor}
 $\asdim(\cH(\Gamma))\leq \lasdim(\Gamma) +1$.
\end{proposition}

\begin{proof}
 We will use the Hurewicz theorem for asymptotic dimension.
\begin{theorem}[{\cite[Theorem 1]{BeDr-Hurewicz}}]
 Let $f:X\to Y$ be a Lipschitz map between geodesic metric spaces, and suppose that for each $R$ the family $\{f^{-1}(B_R(y))\}_{y\in Y}$ has uniform asymptotic dimension $\leq m$. Then $\asdim(X)\leq \asdim Y+m$.
\end{theorem}

Choose some vertex $x \in \Gamma$, and let $\gamma$ be the geodesic ray in $\cH(\Gamma)$ obtained by concatenating 
the edges between $(x,i)$ and $(x,i+1)$ for all $i\in \N$.
Let $f:\cH(\Gamma)\to \gamma$ be the natural $1$-Lipschitz map given by $f(z,i)=(x,i)$.
Fix some $R>0$. The preimages under $f$ of all balls of radius $R$ in $\gamma$ are quasi-isometric 
(with constants depending on $R$ only) to $(\Gamma,d_n)$ for some $n$, 
where $d_n(x,y)=2\ln (d_\Gamma(x,y)e^{-n}+1)$, see Lemma \ref{distestim}. 
So, we only need to show that the uniform asymptotic dimension of $\{(\Gamma,d_n)\}$ is at most $m$. 
Fix any $R$ large enough. There exists $C=C(\Gamma)$ so that the following holds. Let $m=\lasdim(\Gamma)$. 
For each $n$ there exists a covering $\cU(n)=\cU_1(n)\cup\dots\cup \cU_{m+1}(n)$ of $(\Gamma,d_\Gamma)$ such that 
each $\cU_i(n)$ is $e^{n+R}$-disjoint and $Ce^{n+R}$-bounded. 
So, $\cU_i(n)$ is $(2R-M)$-disjoint and $(2R+M)$-bounded with respect to $d_n$, where $M$ depends on $\Gamma$ only.
\end{proof}

Proposition \ref{asdimhor} has a weak converse.
\begin{proposition}\label{prop-asdimhor-converse}
 $\lasdim(\Gamma)\leq \asdim(\cH(\Gamma))$.
\end{proposition}

\begin{proof}
 Set $m=\asdim(\cH(\Gamma))$. There exists a covering $\cU = \cU_1\cup\dots\cup \cU_{m+1}$ of $\cH(\Gamma)$ 
 such that each $\cU_i$ is $r$-disjoint and $R$-bounded, for some $r,R$ large enough. 
 Up to increasing $R$ and decreasing $r$ we have that for each $n$ there exists a covering 
 $\cU_1(n)\cup\dots\cup \cU_{m+1}(n)$ of $(\Gamma,d_n)$ with the same properties,
 where $d_n(x,y)=2\ln (d_\Gamma(x,y)e^{-n}+1)$ as in the previous proof. 
 So, each $\cU_i(n)$ is $e^n(e^{r/2}-1)$-disjoint and $e^n(e^{R/2}-1)$-bounded in $(\Gamma, d_\Gamma)$. 
 We are done as $$ \frac{e^n(e^{R/2}-1)}{e^n(e^{r/2}-1)}$$ is bounded independently of $n$.
\end{proof}

Now we can bound the asymptotic dimension of the Bowditch space.
\begin{proposition}\label{prop-asdim-bowditch}
	Let $G$ be hyperbolic relative to $H_1,\dots, H_n$, and
	let $m = \max_{i=1, \ldots, n} \lasdim(H_i)$.
	Then
	\[
		\max \{ \asdim(G), m \} \leq \asdim(X(G)) \leq \max \{\asdim(G), m+1\}.
	\]
\end{proposition}
\begin{proof}
 The lower bound by $m$ follows from Proposition~\ref{prop-asdimhor-converse}.
 The lower bound by $\asdim(G)$ follows as $X(G)$ contains $G$ with a proper metric.
 
 For the upper bound, we will use the union theorem for asymptotic dimension.
\begin{theorem}[{\cite[Theorem 1]{BeDr-unionthm}}]
 Let $Y$ be a geodesic metric space and suppose that $Y=\bigcup_{i\in\N} A_i$. Also, suppose that $\{A_i\}$ has uniform asymptotic dimension $\leq n$ and for each $R$ there exists $Y_R\subseteq Y$ so that $\asdim Y_R\leq n$ and $\{A_i\backslash Y_R\}$ is $R$-disjoint. Then $\asdim Y\leq n$.
\end{theorem}

We apply the theorem with $\{A_i\}$ the family of horoballs of $X(G)$ and 
for each $R$, take $Y_R$ to be a suitable neighborhood of an orbit of $G$.
By Proposition \ref{asdimhor}, $\asdim(A_i) \leq m+1$, and the constants are uniform as there
are only finitely many different isometry types.
Finally, $\asdim(Y_R) = \asdim(G)$ because the action of $G$ on $X(G)$ is proper.
\end{proof}

\subsection{Linearly controlled metric dimension estimate}
\begin{proposition}\label{prop-cdim-asdim}
	Suppose $X$ is a Gromov hyperbolic geodesic metric space.  Then $\ldim(\bdry X) \leq \asdim(X)$.
\end{proposition}
\begin{proof}
	We fix the notation of Subsection \ref{visual-metric} regarding visual metrics,
	and write $d$ for the metric on $X$.
	
	Each $x \in \bdry X$ is the limit of some $(1,20\delta)$-quasigeodesic $\gamma$ \cite[Remark 2.16]{KB-02-boundaries-survey}.
	Given $R > 0$, define the projection $\pi_R : \bdry X \ra X$ by $\pi_R(x) = \gamma(R)$.
	This is well defined up to an error of $C_1=C_1(\delta)$.
	By considering a quasi-geodesic triangle between $a, b$ and $0$, observe that
	there exists $C_2=C_2(C_0, C_1, \delta)$ so that 
	if $d(\pi_R(a), \pi_R(b)) \geq s > 2C_1$, then  
	$\rho(a,b) \geq \frac1{C_2} e^{-\eps(R-s/2)}$.  
	Similarly, if $d(\pi_R(a), \pi_R(b)) \leq t$ then $\rho(a,b) \leq C_2 e^{-\eps(R-t/2)}$.
	
	If $\asdim(X) \leq n$, then, given $s = 3C_1$, there exists $t < \infty$ and a cover 
	$\cU = \bigcup_0^n \cU_i$ so that each $U \in \cU$ has diameter at most $t$, and $\cU_i$ is $s$-disjoint.
	
	Suppose some small $r>0$ is given.  Let $R = -\frac{1}{\eps} \ln r$.
	For $U \in \cU_i \subset \cU$, let $\hat{U} \subset \bdry X$ be the set of points $z$ so that there exists a 
	$(1,20\delta)$-quasigeodesic $\gamma$ from $0$ to $z$ with $\gamma(R) \in U$.
	Let $\hat{\cU_i} = \bigcup_{U \in \cU_i} \hat{U}$, 
	and $\hat{\cU} = \bigcup_{i=1}^n \hat{\cU_i}$.
	
	By the above estimates, $\hat{\cU_i}$ is $(e^{\eps s/2} r/C_2)$-disjoint, and
	$\hat{\cU}$ is $C_2 r e^{\eps t/2}$-bounded.
	Since the ratio of these distances is bounded by $C_2^2 e^{\eps(t-s)/2}$, we have $\ldim(\bdry X) \leq n$.	
\end{proof}

\section{Embedding in a product of quasi-trees and the Bowditch space}\label{sec-embed-prod-trees-and-bow}
The aim of this section is to find an embedding of a given relatively hyperbolic group into 
a product of trees ``stabilized'' by the Bowditch space or the coned-off graph.
\begin{theorem}\label{prodqtrees}
  Let $G$ be hyperbolic relative to $H_1,\dots, H_n$.
  Suppose there exists $k$ so that
  each $H_i$ admits quasi-isometric embedding into the product of $k$ trees. 
  Then $G$ admits a quasi-isometric embedding into the product of $k$ trees and either $X(G)$ or the coned-off graph $\hat{G}$.
\end{theorem}

We prove this theorem in subsection~\ref{ssec-prodqtrees}.

\subsection{Quasi-trees of spaces}\label{ssec-dist}

To prove Theorem~\ref{prodqtrees} we use a result by Bestvina, Bromberg and Fujiwara described below.

Let $\bf{Y}$ be a set and for each $Y\in\bf{Y}$ let $\cC(Y)$ be a geodesic metric space. 
For each $Y$ let $\pi_Y:{\bf Y}\backslash\{Y\}\to\cP(\cC(Y))$ be a function (where $\cP(\cC(Y))$ is the collection of all subsets of $\cC(Y)$). Define
$$d^{\pi}_Y(X,Z)=\diam\{\pi_Y(X)\cup\pi_Y(Z)\}.$$
Using the enumeration in \cite[Sections 3.1, 2.1]{BBF}, consider the following Axioms. There exists $\xi < \infty$ so that:
\begin{itemize}
	\item[(0)] $\diam(\pi_Y(X))<\xi$ for all distinct $X, Y \in {\bf Y}$,
	\item[(3)] for all distinct $X, Y, Z \in {\bf Y}$ we have $\min\{d^{\pi}_Y(X,Z),d^{\pi}_Z(X,Y)\}\leq \xi$,
	\item[(4)] $\{Y:d^\pi_Y(X,Z)\geq \xi\}$ is a finite set for each $X,Z\in\bf{Y}$.
\end{itemize}

Using the functions $d_Y^\pi$, the authors of \cite{BBF} define a certain complex
$\cP_K(\{(\cC(Y),\pi_Y)\}_{Y\in\bf{Y}})$ with vertex set ${\bf Y}$, which they denote as $\cP_K({\bf Y})$,
and turns out to be a quasi-tree \cite[Section 2.3]{BBF}.

Let $\cC(\{(\cC(Y),\pi_Y)\}_{Y\in\bf{Y}})$ be the path metric space
consisting of the union of all $\cC(Y)$'s and edges of length 1 connecting all points in $\pi_X(Z)$ to 
all points in $\pi_Z(X)$ whenever $X,Z$ are connected by an edge in the complex
$\cP_K(\{(\cC(Y),\pi_Y)\}_{Y\in\bf{Y}})$.
More details of this construction are found in \cite[Section 3.1]{BBF}.

The only result from \cite{BBF} that we need is the following.

\begin{theorem}[{\cite[Theorem 3.10]{BBF}}]\label{bbf}
 If $\{(\cC(Y),\pi_Y)\}_{Y\in\bf{Y}}$ satisfies Axioms $(0), (3)$ and $(4)$ and each $\cC(Y)$ is a tree then $\cC(\{(\cC(Y),\pi_Y)\}_{Y\in\bf{Y}})$ is a quasi-tree (i.e. it is quasi-isometric to a tree).
\end{theorem}

To apply this theorem in the context of Theorem~\ref{prodqtrees}, we fix a Cayley graph $\Gamma(G,S)$ 
with word metric $d$ of the relatively hyperbolic group $G$ as in Definition~\ref{def-rel-hyp}.
Now let $\bf{Y}$ be the collection of all left cosets of peripheral subgroups in $G$.
For each such left coset $Y = g H_i$ we denote by $\cC(Y)$ the subgraph $g \Gamma(H_i, S \cap H_i)$ in $\Gamma(G,S)$,
which contains $Y$ and is a copy of the Cayley graph of $H_i$.
Let $\pi_Y$ be a closest point projection on $\cC(Y)$ for each $Y\in\bf{Y}$.
Namely, for each $x \in X$, let $\pi_Y(x)$ be the subset of $\cC(Y)$ minimizing the function $d(x, \cdot)$,
and let $\pi_Y(X) = \bigcup_{x \in X} \pi_Y(x)$.

As peripheral subgroups are undistorted in $G$ \cite[Lemma 4.15]{DSp-05-asymp-cones} we can, for our purposes, 
identify $\cC(Y)$ with the corresponding subset of (the Cayley graph of) $G$.
\begin{lemma}
	The collection $\{(\cC(Y),\pi_Y)\}_{Y\in\bf{Y}}$ satisfies Axioms $(0), (3)$ and $(4)$.
\end{lemma}
\begin{proof}
	For some constant $C$, the following holds \cite{projrelhyp}:
\begin{itemize}
 \item $\diam(\pi_Y(Y'))\leq C$ whenever $Y\neq Y'$, and
 \item if for some $x,y \in G$ we have $d(\pi_Y(x),\pi_Y(y))\geq C$ then any geodesic from $x$ to $y$ in the Cayley graph of $G$ intersects
 both of the balls of radius $C$ around $\pi_Y(x)$ and $\pi_Y(y)$.
\end{itemize}

The first property clearly implies Axiom $(0)$. The fact that the second property implies Axiom $(3)$ can be considered folklore, see \cite[Lemma 2.5]{contr}. Axiom $(4)$ follows from the fact that the right hand side of the distance formula \eqref{eqn-distform} is finite.
\end{proof}

\subsection{The proof of Theorem~\ref{prodqtrees}}\label{ssec-prodqtrees}
Let us go back to the general setting of Theorem \ref{bbf} for a moment, with 
$\{ \cC(Y), \pi_Y \}_{Y \in {\bf Y}}$ as in the statement of that theorem.
Suppose that we are also given another collection of geodesic metric spaces $\{ \cC'(Y) \}_{Y \in {\bf Y}}$,
indexed by the same set ${\bf Y}$, and coarsely Lipschitz maps 
$f_Y:\cC(Y)\to \cC'(Y)$, for $Y \in {\bf Y}$, with uniform constants. 
Observe that if axioms $(0), (3)$ and $(4)$ hold for $\{(\cC(Y),\pi_Y)\}_{Y\in\bf{Y}}$, then 
they also hold for $\{(\cC'(Y),f_Y\circ\pi_Y)\}_{Y\in\bf{Y}}$.

Now suppose that each $H_i$ admits a quasi-isometric embedding into the product of $k$ trees. 
We then have coarsely Lipschitz maps $f_{i,Y}$, $i=1,\dots,k$ from $\cC(Y)$ to some tree $T_{i,Y}$. 
According to the observation we just made and Theorem \ref{bbf}, for each $i$ we have a map 
$p_i: G\to {\bf T}_i=\cC(\{(T_{i,Y},f_{i,Y}\circ \pi_Y)\})_{Y\in {\bf Y}}$, and ${\bf T}_i$ is a quasi-tree. 
The last step in the proof of Theorem \ref{prodqtrees} is the following proposition.
\begin{proposition}
 The map $f=\prod p_i\times c: G\to \prod {\bf T}_i\times \hat{G}$ is a quasi-isometric embedding, 
 where $c:G\to \hat{G}$ is the inclusion.
\end{proposition}
\begin{proof}
There are distance formulas for the ${\bf T}_i$'s \cite[Lemma 3.3, Corollary 3.12]{BBF} which summed up give,
for $L,\lambda,\mu$ large enough,
$$d\left(\prod p_i(x),\prod p_i(y)\right)\approx_{\lambda,\mu} \sum_{Y\in\bf{Y}} \big\{\big\{d(\pi_Y(x),\pi_Y(y))\big\}\big\}_L.$$
Comparing this with the distance formula for relatively hyperbolic groups (Theorem~\ref{distform})
immediately gives the required estimates.
\end{proof}

The map $c$ factors through Lipschitz maps $G\to X(G)\to \hat{G}$, so the proposition implies both versions of Theorem \ref{prodqtrees}.

\section{Three dimensional manifold groups}\label{sec-three-mfld}

In this section we prove the following theorem.
\begin{theorem}\label{thm-three-mfld}
	Let $G = \pi_1(M)$, where $M$ is a compact, orientable $3$-manifold whose (possibly empty) boundary is a union of tori.
	Then $\ecodim(G) < \infty$ if and only if no manifold in the prime decomposition of $M$ has $\mathrm{Nil}$ geometry;
	in this case, $\ecodim(G) \leq 8$.
	
	In any case, $\lasdim(G) \leq 8$.
\end{theorem}

The second part of the following proposition is not needed for our purposes, but we think it is of independent interest.

\begin{proposition}\label{nonclosed}
 \begin{enumerate}
  \item If $M$ is a graph manifold with non-empty boundary then $\lasdim(\pi_1(M))=2$.
  \item If $M$ is a Haken orientable $3$-manifold with non-empty boundary then $\asdim(\pi_1(M))\leq 2$.
\end{enumerate}
\end{proposition}

Recall that a compact orientable $3$-manifold $M$ is \emph{Haken} if it is irreducible and it contains a $\pi_1$-injective embedded surface $S$. Cutting $M$ along $S$ gives another Haken manifold, and moreover the new $\pi_1$-injective surface can be required to have non-empty boundary. Repeating the cutting procedure finitely many times yields a disjoint union of balls. Such sequence of cuts is called \emph{Haken hierarchy}.

\begin{proof}
 $1)$ The lower bound follows from the existence of undistorted copies of $\mathbb{Z}^2$ in $\pi_1(M)$. Up to passing to a finite-sheeted cover of $M$, we can assume that $M$ fibers over $S^1$ by \cite{WaYu-graphfiber}, so that we have a short exact sequence
$$1\to F\to \pi_1(M)\to \mathbb{Z}\to 1,$$
where $F$ is a free group.
 What is more, the content of \cite[Theorem 0.7]{WaYu-graphfiber} is that a fiber intersects all Seifert components and it splits them open into products of a surface with an interval. In particular, the monodromy of the fiber bundle is a product of Dehn twists along disjoint simple closed curves (the connected components of the intersection of the fiber with the boundary of the Seifert components), and so $F$ is undistorted in $\pi_1(M)$. The result then follows from a direct application of \cite[Theorem 0.2]{Bretal-Hurewicz-AN}.
\par
$2)$ We proceed by induction on the length of a Haken hierarchy for $M$ so that all surfaces involved have non-empty boundary. If $M$ is a ball, then the result trivially holds. Otherwise, $\pi_1(M)$ can be written as an amalgamated product $\pi_1(N_1)*_{\pi_1(S)}\pi_1(N_2)$ or an HNN extension $\pi_1(N_1)*_{\pi_1(S)}$, where $N_i$ is a Haken manifold admitting a strictly shorter Haken hierarchy (of the type described above) than $M$ and $S$ is a compact surface with non-empty boundary. We can then use the induction step and the fact that, when $A,B,C$ are finitely generated groups, $\asdim(A*_C B)\leq \max\{\asdim A,\asdim B,\asdim C+1\}$ and $\asdim(A*_C)\leq \max\{\asdim A,\asdim C+1\}$ (see \cite{Dr-amalg-form} and Theorem \ref{asdimhnn}).
\end{proof}

\begin{remark}
 The second part of the proposition (in the case of toric boundary) also follows from recent results about virtual specialness of fundamental groups of $3$-manifolds \cite{Wi-qconv-hier,Liu-virtcubgraph,PW-mixedvirtspec} in view of a criterion for virtual fibering discovered by Agol \cite{Ag-critvirtfibr}.
\end{remark}

We now show the following easy lemma and then prove the theorem.

\begin{lemma}
 $\ecodim G_1 *G_2=\max\{\ecodim G_1,\ecodim G_2\}$.
\end{lemma}

\begin{proof}
 The inequality $\geq$ follows from the fact that $G_1,G_2$ are undistorted in $G=G_1*G_2$. Suppose that, for $i=1,2$, $f_i:G_i\to \prod_{j=1}^n T^i_j$ are quasi-isometric embeddings. We have to show that $G$ embeds in the product of $n$ trees as well. Denote by $T$ the Bass-Serre tree of $G$. For each vertex $v$ of $T$ denote by $T_k^v$ a copy of $T^{i(v)}_k$, where $i(v)$ equals $i$ if the stabilizer of $v$ is conjugate to $G_i$.
When $e$ is an edge of $T$ with endpoints $v_1$, $v_2$, we let $p_e$ be the only element of $G$ in the intersection of the left cosets of $G_1, G_2$ corresponding to $v_1,v_2$. Finally, we let $T_k$ be the tree obtained from $\bigcup T_k^v$ by adding an edge of length, say, $1$ connecting $f_{1,k}(p_e)$ to $f_{2,k}(p_e)$ for each edge $e$ of $T$, where $f_{i,k}$ is the $k$-th component of $f_i$ and we identify $G_i$ with its left coset corresponding to an endpoint of $e$.
\par
There is a natural map $f:G\to\prod_{i=1}^n T_k$, which (up to an error bounded by $1$) restricts to $f_i$ on every left coset of $G_i$. It is readily checked that this map is a quasi-isometric embedding (using more sophisticated technology than needed, one can use the distance formula and observe that, as we added edges of length $1$ connecting $T_{v_1}$ to $T_{v_2}$ when $v_1$ is adjacent to $v_2$, we have $d(f(x),f(y))\geq d_T(x,y)$ and $d_T(x,y)$ is approximately $d_{\hat{G}}(x,y)$).
\end{proof}

\begin{proof}[Proof of Theorem \ref{thm-three-mfld}]
	Let $G_i$, for $i=1, \ldots, n$, be the fundamental groups of the prime factors $M_i$ of $M$.
	Then $G$ is the free product of the $G_i$, so that $\ecodim(G)=\max\{\ecodim(G_i)\}$ in view of the previous lemma, and also $\lasdim(G) = \max \{ \lasdim(G_i) \}$ by \cite{BrHi-09-AN-tree-graded}. In particular, we can just study the case when $M$ is prime.
\par
We report below a list of previously known cases. Except for the last two cases, the first column indicates the geometry (of the universal cover) of $M$.
The values in the table are justified below.
	\begin{center}
	\begin{tabular}{c|c|c}
		$M$ & $\lasdim(\pi_1(M))$ & $\ecodim(\pi_1(M))$ \\
		\hline
		$S^3$			& $0$	& $0$ \\
		$\R^3$			& $3$	& $3$ \\
		$\HH^3$, closed			& $3$	& $3$ \\
		$S^2 \times \R$	& $1$	& $1$ \\
		$\HH^2 \times\R$, closed & $3$	& $3$ \\
		$\HH^2 \times\R$, non-closed & $2$	& $2$ \\
		$\widetilde{SL_2\R}$ & $3$	& $3$ \\
		$\mathrm{Nil}$	& $3$	& $\infty$ \\
		$\mathrm{Sol}$	& $3$	& $3$ or $4$\\
graph manifold, closed & $3$ & $3$\\
graph manifold, non-closed & $2$ & $2$ or $3$
	\end{tabular}
	\end{center}
	The bounds on $S^3$, $S^2 \times \R$ and $\R^3$ are trivial.
	The calculation of $\lasdim(\HH^n) = \ecodim(\HH^n) = n$ is found in 
	\cite{BuSh-05-hyp-tree-bound-above,BuSh-07-hyp-tree-bound-below}.
	This also gives $\ecodim(\HH^2 \times \R) \leq 3$, and $\lasdim(\HH^2 \times \R) = 3$ 
	\cite[Theorem 4.3]{DrSm-07-an-asdim-prod-r}.
	The spaces $\HH^2 \times\R$ and $\widetilde{SL_2\R}$ are quasi-isometric. If $M$ is not closed and has geometry $\HH^2 \times\R$ then $\pi_1(M)$ is virtually the product of a free group and $\mathbb{Z}$. 
	
	The discrete Heisenberg group $H$ is quasi-isometric to $\Nil$.
	The result $\asdim(H)=\lasdim(H)=3$ has been obtained by several people, for example see \cite{DyHi-08-ascone-an-dim}.
	On the other hand, $\ecodim(H)=\infty$ as $H$ does not admit a quasi-isometric embedding into the product 
	of finitely many metric trees, or indeed any CAT(0) space \cite{Pauls-01-qi-heis}.
	Thus the proof in the ``only if'' direction is complete.

	For $\lasdim(\Sol)=3$ see \cite{HiPe-09-andim-lie}.
	As $\Sol$ quasi-isometrically embeds in $\HH^2 \times \HH^2$ (see, for example, \cite[Section 9]{deCor-08-dim-as-cone-lie}),
	$\ecodim(\Sol) \leq 4$.
	
	If $M$ is a graph manifold then $\lasdim(\pi_1(M))\leq\ecodim(\pi_1(M)) \leq 3$ by \cite{HuSi-11-andim-graph-manifold}, and it is observed in the same paper that the equalities hold if $M$ is closed. We handled the non-closed case for $\lasdim$ in Proposition \ref{nonclosed}.

There are only two cases left. First, when $M$ is a finite volume, non-closed hyperbolic manifold
it is well-known that $\pi_1(M)$ is hyperbolic relative to virtually $\mathbb{Z}^2$ 
subgroups \cite{Farb-98-rel-hyp}.  (Notice also that $\asdim(\pi_1(M))\leq 3$ as $\pi_1(M)$ admits a coarse embedding in $\HH^3$.)
As $X(\pi_1(M))$ is quasi-isometric to $\HH^3$,
by Theorem~\ref{prodqtrees},
\[
	2 \leq \lasdim(\pi_1(M)) \leq \ecodim(\pi_1(M)) \leq 5.
\]

The second case is when $M$ 
is non-geometric and its geometric decomposition contains at least one hyperbolic component. 
In this case $\pi_1(M)$ is hyperbolic relative to virtually $\mathbb{Z}^2$ and graph manifold groups,
as a consequence of the combination theorem \cite[Theorem 0.1]{Dah-03-combination} 
and aforementioned fact that the hyperbolic component groups are hyperbolic relative to virtually $\Z^2$ groups
(a full statement is given in \cite[Theorem 9.2]{AFW-12-3manifoldgroups}).
Therefore, by the graph manifold groups bound of \cite{HuSi-11-andim-graph-manifold}, every peripheral group
has $\ecodim$ at most $m=3$.
Note too that $\asdim(\pi_1(M))\leq 4$ holds by \cite{BeDr-04-asdim-trees}.
Thus by Theorem~\ref{thm-main1},
\begin{align*}
	\ecodim(\pi_1(M) & \leq \max\{\asdim(\pi_1(M)),m+1\}+m+1 \\ & = \max\{ 4, 3+1\} + 3+ 1 = 8.
\end{align*}
Therefore, we have
\[
	2 \leq \lasdim(\pi_1(M)) \leq \ecodim(\pi_1(M)) \leq 8,
\]
when $M$ is non-closed, and
\[
	3 \leq \lasdim(\pi_1(M)) \leq \ecodim(\pi_1(M)) \leq 8,
\]
when $M$ is closed.  The lower bound in the last case follows from the fact that the asymptotic dimension is 
greater or equal to the virtual co-homological dimension for groups of type FP \cite{Dr-cohom-asym}.
\end{proof}

\appendix
\section{Asymptotic dimension of HNN extensions}\label{sec-hnn}

The proof of the following theorem follows almost verbatim the arguments in \cite{Dr-amalg-form}.

\begin{theorem}
\label{asdimhnn}
 Let $A,C$ be finitely generated groups. Then
$$\asdim(A*_C)\leq\max\{\asdim(A),\asdim(C)+1\}.$$
\end{theorem}

We require the following notation.  For $X$ a metric space, we will say that $(r,d)$-$\dim X\leq n$ if there exists a $d$-bounded cover of $X$ with Lebesgue number at least $r$ and $0$-multiplicity at most $n+1$. We will say that a family $\{X_i\}$ of metric spaces satisfy $\asdim X_i\leq n$ uniformly if for every $r>0$ there exists $d$ so that $(r,d)$-$\dim X_i\leq n$. A partition of the metric space $X$ is a presentation of $X$ as a union of subspaces with pairwise disjoint interiors.  The proof of Theorem~\ref{asdimhnn} uses the criterion below.

\begin{theorem}[{\cite[Partition Theorem]{Dr-amalg-form}}]\label{thm-partition}
 Let $X$ be a geodesic metric space. Suppose that for every $R>0$ there exists $d>0$ and a partition $X=\bigcup_{i\in\N} W_i$ where $\asdim W_i\leq n$ uniformly and with the property that $(R,d)$-$\dim(\bigcup \partial W_i)\leq n-1$ (with the restriction of the metric of $X$). Then $\asdim (X)\leq n$.
\end{theorem}

We now do some preliminary work.
\par
 Fix a generating system $S_A$ of $A$ and let $t$ be the stable letter of the HNN extension.  Let $d$ be the associated word metric on $G=A *_C$.
 The group $G$ acts on its Bass-Serre tree whose vertices are left cosets of $A$ in $G$ and whose edges are labeled by left cosets of $C$. The endpoints of the edge $gC$ are $gA$ and $gtA$. Denote by $K$ the graph dual to the Bass-Serre tree. We will denote the simplicial metric in $K$ by $|\cdot,\cdot|$, and we let $|u|=|u,C|$. Notice that for each pair of vertices in $K$ there exists a unique geodesic connecting them. Let $\pi:G\to K$ be the map $g\mapsto gC$.
\begin{remark}
\label{1lip}
 $\pi$ extends to a simplicial map of the Cayley graph of $G$. In particular, $\pi$ is $1$-Lipschitz.
\end{remark}

In fact, let $s\in S_A$. Then for each $g\in G$ we have $gsA=gA$, so that the edges $gC$ and $gsC$ of the Bass-Serre tree share the endpoint $gA$, which is exactly saying that the vertices $gC, gsC$ of $K$ are connected by an edge. Similarly, the edges $gC$ and $gtC$ of the Bass-Serre tree share the endpoint $gtA$.
\par
We divide $K$ into two parts $K_0,K_1$ intersecting only at the base vertex $C$, where $K_1$ contains the edges corresponding to $tA$.  Let $\overline{d}$ be the graph metric on $K$, and denote by $B_r^1$ the closed $r$-ball in $K_1$ centered at $C$, where $r$ will always denote an integer. We will write $v\leq u$, where $u,v\in K^{(0)}$, if $v$ lies on the geodesic segment $[C,u]$ (notice that this is a partial order). For $u\in K^{(0)}$ with $u\neq C$ and $r>0$ denote
$$K^u=\{v\in K^{(0)}|\ v\geq u\},$$
and
$$B^u_r=\{v\in K^u|\ |v|\leq |u|+r\}.$$
Notice that if $u=gC \in K_1$ then $B^u_r=gB^1_r$ and $K^u=gK_1$. In particular, $\pi^{-1}(B^u_r)=g\pi^{-1}(B^1_r)$ and $\pi^{-1}(K^u)=g\pi^{-1}(K_1)$.
\par
We say that $F\subseteq G$ separates $H_1,H_2\subseteq G$ if all paths in the Cayley graph connecting $H_1$ to $H_2$ intersect $F$. Set $D_R=\{x\in G|\ d(x,C)=R\}\cap \pi^{-1}(K_1)$, for $R\in \N$. For $u=gC\in K_1$ denote $D^u_R=g D_R$ and notice that $\pi(D^u_R)\subseteq B^u_R$.
\begin{lemma}
\label{tech1}
 Let $u \in K^{(0)}_0, v\in K^{(0)}$ be so that either $v<u$ or $v$ is incomparable with $u$. Then $D^u_R$ separates $\pi^{-1}(v)$ and $\pi^{-1}(u')$ whenever $u<u'$ and $|u'|-|u|>R$.
\end{lemma}

\begin{proof}
 We will show that $D_R$ separates $\pi^{-1}(K_0)$ and $\pi^{-1}(u')$ if $|u'|>R$ and $u'\in K_1$, using the action of $G$ then yields the required statement. As $\pi$ is $1$-Lipschitz we have $d(C,\pi^{-1}(u'))>R$ so that $D_R$ separates $C$ and $\pi^{-1}(u')$. To complete the proof notice that $C$ separates $\pi^{-1}(K_0)$ and $\pi^{-1}(K_1)$ (as their images through $\pi$ are separated by the vertex labeled $C$).
\end{proof}

\begin{lemma}
\label{tech2}
 Suppose $R\leq r/4$. Then $d(gD_R,g'D_R)\geq 2R$ whenever $g,g'\in G$ with $gC, g'C \in K_1$, $|gC|,|g'C|\in r\N$ and $gC\neq g'C$.
\end{lemma}

\begin{proof}
 Set $u=gC, u'=g'C$. Suppose first that $|u|\neq|u'|$. As $\pi(gD_R)\subseteq B^u_R$ and $\overline{d}(B^u_R,B^{u'}_R)\geq r-R\geq 3R$, and in view of the fact that $\pi$ is $1$-Lipschitz, we get $d(gD_R,g'D_R)\geq 3R$.
\par
Suppose instead $|u|=|u'|$ and pick $x\in g D_R, y\in g'D_R$. Every path in $K$ between $\pi(x)$ and $\pi(y)$ passes through $u$ and $u'$. This applies in particular to the projection of a geodesic $\gamma$ in the Cayley graph of $G$, so that $\gamma$ intersects $gC$ and $g'C$. Since $d(x,gC),d(y,g'C)=R$, we get $d(x,y)\geq 2R$.
\end{proof}

\begin{lemma}
 For $m\in \N$, let $(At)^m= \bigcup A t^{\eps_1} \cdots A t^{\eps_m} \subseteq A *_C$, where each $\eps_i$ equals $0$ or $1$.  Then $\asdim(At)^m\leq\max\{\asdim A,\asdim C+1\}$ for every $m$.
\end{lemma}

\begin{proof}
$(At)^m$ admits a coarse embedding in $\pi_1(\cG)$, where $\cG$ is a graph of groups so that all vertex groups are isomorphic to $A$, all edge groups are isomorphic to $C$ and the underlying graph is a tree. Repeated applications of \cite[Theorem 2.1]{Dr-amalg-form} give $\asdim \pi_1(\cG)\leq \max\{\asdim A,\asdim C+1\}$, and hence the same holds for $(At)^m$.
\end{proof}

We are now ready for the proof of the theorem.

\emph{Proof of Theorem \ref{asdimhnn}.} Set $n=\max\{\asdim(A),\asdim(C)+1\}$. Once we show $\asdim(\pi^{-1}(K_0)),\asdim(\pi^{-1}(K_1)) \leq n$, the conclusion follows from the Finite Union Theorem \cite{BeDr-unionthm}. We will show the latter, using the Partition Theorem \ref{thm-partition}. Fix $R>0$ and take $r>4R$. By Lemma \ref{tech1} we can write $G=X_+\cup X_{-}$ with $X_+\cap X_{-}=D_R$ so that $X_+\subseteq \pi^{-1}(K_1)$, and $\pi^{-1}(K_0)\subseteq X_{-}$, and $D_R$ separates $X_+\backslash D_R$ and $X_{-}\backslash D_R$. For each $u\in K_1^{(0)}$ fix $g_u$ so that $u=g_u C$. Set $X^u_{\pm}=g_u(X_{\pm})$, $V_r=X_+\cap \left(\bigcap_{|u|=r} X^u_{-}\right)$ and $V^u_r=g_u(V_r)$. It is readily seen that $\pi(V_r)\subseteq B_{r+R}^1$ and that
$$\pi^{-1}(K_1)=\bigcup_{|u|\in r\N^+} V^u_r\cup N^+_R(C),$$
where $N^+_R(C)=N_R(C)\cap\pi^{-1}(K_1)$. If $|u|,|w| \in r\N^+$, $V^u_r\cap V^w_r\neq \emptyset$ and $u\neq w$ then either $u<w$ and $|w|=|u|+r$ or, vice versa, $w<u$ and $|u|=|w|+r$. Also, if $V^u_r\cap V^w_r\neq \emptyset$ and $u<w$ then $V^u_r\cap V^w_r=D^w_R$ for $D^w_R=g_wD_R$. Putting these facts together we get
$$Z=\bigcup_{|u|\in r\N^+} \partial V^u_r=\bigcup_{|u|\in r\N^+} D^u_R.$$
As $D_R$ is coarsely equivalent to $C$, there is, for some $d>0$, an $(R,d)$-cover $\cU$ of $D_R$ (i.e. the Lebesgue number of $\cU$ is at least $R$ and $\cU$ is $d$-bounded)  with $0$-multiplicity at most $n$. By Lemma \ref{tech2}, $\bigcup_{|u|\in r\N^+} g_u\cU$ is an $(R,d)$-cover of $Z$, and this witnesses the fact that $(R,d)$-$\dim(Z)\leq n-1$. Finally, notice that $\pi^{-1}(B_s^1)\subseteq (At)^{s+1}$ so that $\asdim B_s^1\leq n$, for each $s\in\N$. In particular, $\asdim \pi^{-1}(B_{r+R}^1)\leq n$ and thus $\asdim \pi^{-1}(V^u_r)\leq n$ uniformly. Finally, $\asdim N^+_R(C)\leq n-1\leq n$, so that the Partition Theorem applies.\qed

\bibliographystyle{alpha}
\bibliography{biblio}

\end{document}